\documentclass[preprint]{elsarticle}
\usepackage{amssymb,url,amsmath,mathtools}
\usepackage[english]{babel}

\newtheorem{theorem}{Theorem}[section] % 1st argument is your name for it
\newtheorem{lemma}[theorem]{Lemma}     % 2nd argument is what is printed
\newtheorem{corollary}[theorem]{Corollary}

\newtheorem{conj}[theorem]{Conjecture}
\def\T{{\mathcal T}}
\def\C{{\mathbb C}}
\def\F{{\mathbb F}}
\def\Q{{\mathbb Q}}
\def\Z{{\mathbb Z}}
\newcommand\myatop[2]{\genfrac{}{}{0pt}{}{#1}{#2}}

\newproof{pf}{Proof}

\begin{document}

\begin{frontmatter}
\title{The Eigenvalues of the Graphs $D(4,q)$}

\author[uw]{G.\ Eric Moorhouse\corref{cor1}}
\ead{moorhous@uwyo.edu}

\author[ud]{Shuying Sun}
\ead{shuying@udel.edu}

\author[uw]{Jason Williford}
\ead{jwillif1@uwyo.edu}

\cortext[cor1]{Corresponding author}
\address[uw]{Department of Mathematics, University of Wyoming, Laramie WY 82071 USA}
\address[ud]{Department of Mathematical Sciences, University of Delaware, Newark DE 19716 USA}

\begin{abstract}
The graphs $D(k,q)$ have connected components $CD(k,q)$ giving the best known bounds on extremal
problems with {\em forbidden\/} even
cycles, and are denser than the well-known graphs of Lubotzky, Phillips, Sarnak \cite{LPS}
and Margulis~\cite{M1,M2}. Despite this,
little about the spectrum and expansion properties of these graphs is known. In this paper we find the
spectrum for $k=4$, the smallest open case. For each prime power $q$, the graph $D(4,q)$ is
$q$-regular graph on $2q^4$ vertices, all of whose eigenvalues other than $\pm q$ are bounded in
absolute value by $2\sqrt{q}$. Accordingly, these graphs are good expanders,
in fact very close to Ramanujan.
\end{abstract}
\end{frontmatter}

\section{Introduction}\label{introduction}

Let $\Gamma$ be a graph with vertex set $V$. (All our graphs are undirected and have no loops or
multiple edges. See e.g.~\cite{GR,BH} for standard terminology and theory of graphs.) Given a set of
vertices $F\subset V$, we define $\partial F$ to be the
set of vertices in $V\smallsetminus F$ which are adjacent to some vertex of $F$. The
{\em isoperimetric constant\/} of $\Gamma$ is defined to be
\[h(\Gamma)=\min\left\{\frac{|\partial F|}{|F|} : F \subset V \text{ and } |F|\leqslant\frac{|V|}{2}\right\}.\]
An infinite family of $d$-regular graphs whose isoperimetric constants are uniformly bounded away from 0 is
an {\em expander family\/}. The best known general bounds on $h(\Gamma)$ are expressed in terms of
the {\em spectrum\/} of $\Gamma$, i.e.\ the multiset of eigenvalues of its adjacency matrix.
In particular, if $\Gamma$ is $q$-regular with
second-largest eigenvalue $\lambda_2(\Gamma)<\lambda_1=q$, then
\[{\textstyle{\frac12}}\bigl(q-\lambda_2(\Gamma)\bigr)\leqslant h(\Gamma)\leqslant\sqrt{2q(q-\lambda_2(\Gamma))}\,;\]
see e.g.~\cite[Prop.1.84]{KS}. (The second-largest eigenvalue is denoted differently in some
sources, including \cite{KS}.) Thus to certify an infinite family of $q$-regular graphs as an expander
family, we require a uniform lower bound on the {\em spectral gap\/} $q-\lambda_2(\Gamma)$.
A $q$-regular connected graph $\Gamma$ is {\em Ramanujan\/} if
$\lambda_2(\Gamma)\leqslant2\sqrt{q-1}$; by the Alon-Boppana Theorem
 (see e.g.~\cite[Ch.3]{KS}) this bound is asymptotically best possible for any infinite family of $q$-regular graphs.

In searching for good families of explicitly defined graphs with good expansion, a particularly promising
infinite family of graphs is the sequence
\[\cdots\to D(5,q)\to D(4,q)\to D(3,q)\to D(2,q)\]
defined by Lazebnik and Ustimenko~\cite{LU} for each prime power $q$. Each graph $D(k,q)$
in this sequence is bipartite $q$-regular on $2q^k$ vertices having girth${}\geqslant k+4$
(or $k+5$, when $k$ is even); and each connecting
map `$\to$' is a graph-theoretic cover (see \cite[Sec.3B]{LW}). The graphs $D(k,q)$ are connected for
$k\leqslant5$ and $q$ odd; see~\cite{LUW}. The covering property ensures that the girth of $D(k,q)$ is
weakly increasing as $k\to\infty$, and the spectrum of $D(k,q)$ embeds in that of $D(k+1,q)$; see~\cite[Sec.3C]{LW}.

The graphs $CD(k,q)$ are important in the study of Tur\'an type problems on even cycles, giving better
lower bounds on the maximum number of edges in graphs of girth $g\geqslant6$ than the well-known
Ramanujan graphs of Lubotzky, Phillips and Sarnak~\cite{LPS}. Similarly, the graphs $LD(q,r)$ of
Alon et al.~\cite{ASS}, another expander family with fixed degree, have girth 3 (after removing loops).
By comparison, therefore, one might expect the graphs $CD(k,q)$ to have very good expansion properties.

However, little is known about the eigenvalues of these graphs. In fact, to
date only the spectrum of $D(2,q)$ and $D(3,q)$ are known, their characteristic polynomials being
\[(x^2-q^2)(x^2-q)^{q-1}x^{2q(q-1)}\]
and
\[(x^2-q^2)(x^2-2q)^{q(q-1)^2/2}(x^2-q)^{2q(q-1)}x^{q^3-2q^2+3q-2}\]
respectively; see \cite[Sec.5]{LLW}. In particular, these graphs are Ramanujan. However,
Reichard~\cite{R} and Thomason~\cite{T} independently showed by
computer that the graphs $D(4,q)$ are not Ramanujan for certain $q$, refuting the claim of \cite{U1}; see
also the final note in this paper where we investigate this question more closely.
The same statements apply to $D(k,q)$ for all $k\geqslant4$, since the spectrum
of $D(4,q)$ is embedded in that of $D(k,q)$ for $k\geqslant4$.

It was later claimed in \cite{U2} that the eigenvalues of $D(k,q)$ other than $\pm q$ are bounded by
$2\sqrt{q}$. However, a flaw was later found in the argument, leaving the problem open; see the
Math Review MR2048644 for~\cite{U2}.
To date, we have not found any counterexample to this statement, so we list it as a conjecture.
Following~\cite{LU}, we denote by $CD(k,q)$ a connected component of $D(k,q)$; and we note that
$CD(4,q)=D(4,q)$ whenever $q\notin\{2,4\}$.

\begin{conj}[Ustimenko]\label{conjecture}
For all $(k,q)$, $CD(k,q)$ has second largest eigenvalue less than or equal to $2\sqrt{q}$.
\end{conj}
In this paper we verify Conjecture~\ref{conjecture} for $k=4$:
\begin{theorem}
The second largest eigenvalue of $CD(4,q)$ is less than or equal to $2\sqrt{q}$.
\end{theorem}
This implies that these graphs are very close to Ramanujan. Our proof is given in Section~\ref{q_even}
for even $q$, and in Section~\ref{q_odd} for odd $q$. A more explicit determination of the spectrum
is given in Section~\ref{exact} for prime values $q=p$.

Our approach is similar to \cite{CLL}, in that we first realize the halved (point) graph
of $D(k,q)$ as a Cayley graph of a certain $p$-group~$G$. Unlike the situation for the Wenger graphs
in~\cite{CLL}, or the graphs $D(2,q)$ and $D(3,q)$, our group $G$ is nonabelian whenever $q$ is odd,
thus requiring more extensive use of the representation theory of $G$. Finally, our bounds on eigenvalues
are obtained using Weil's bound for exponential sums over $\F_q$, or over Galois rings of characteristic~3
in the case $q=3^e$.
 
\section{The Graphs $D(4,q)$ and their Point Collinearity Graphs $\Gamma(4,q)$}\label{results}

Throughout, we take $F=\F_q$ where $q$ is a prime power. The graph $D(4,q)$
is bipartite and $q$-regular with $2q^4$ vertices. These include $q^4$ vertices $P=P(p_1,p_2,
p_3,p_4)$ called {\it points\/}, and $q^4$ vertices $L=L(\ell_1,\ell_2,\ell_3,\ell_4)$
called {\it lines\/}, where all coordinates are in $F$; and the point $P$ and line $L$ (with coordinates as
above) are {\it incident\/} iff
\[p_2+\ell_2=p_1\ell_1,\quad p_3+\ell_3=p_1\ell_2\quad\hbox{and}
\quad p_4+\ell_4=p_2\ell_1.\]
Also denote by $\Gamma=\Gamma(4,q)$ the point collinearity graph of $D(4,q)$, i.e.\ the graph whose vertices
are the points of $D(4,q)$, two points being adjacent in $\Gamma$ iff they are distinct but collinear in
$D(4,q)$; see e.g.~\cite[Sec.14.2.2]{BH}. One checks that two vertices $P(p_1,p_2,p_3,p_4),P(p_1',p_2',p_3',p_4')$ are adjacent
in $\Gamma$ (i.e.\ distinct and collinear in $D(4,q)$) iff
\[p_1\neq p_1',\quad(p_1-p_1')(p_4-p_4')=(p_2-p_2')^2\quad\hbox{and}\quad p_3-p_3'=p_2 p_1'-p_1 p_2'.\]
The adjacency matrix of $\Gamma$ has the form $A=B_1 B_1^T-qI_{q^4}$ where $B_1$ is a
$q^4\times q^4$ matrix for which
\[B:=\left[\begin{array}{cc}0&B_1\\ B_1^T&0\end{array}\right]\]
is the adjacency matrix of $D(4,q)$ (with the first $q^4$ rows and columns indexed by points, and the
last $q^4$ rows and columns indexed by lines). Note that $\Gamma$ is
a $q(q-1)$-regular graph on $q^4$ vertices. The spectra of $A$ and $B$
are in direct relationship. Indeed, elementary methods yield the following, which is also implicit in
\cite{CLL,LLW}:

\begin{lemma}\label{AB}
Denote the characteristic polynomial of $A$, the adjacency matrix of $\Gamma(4,q)$, by
$\phi(x)=\det(xI_{q^4}-A)$. Then the characteristic polynomial of $B$, the adjacency matrix of
$D(4,q)$, is\/ $\det(xI_{2q^4}-B)=\phi(x^2-q)$.\qed
\end{lemma}
Equivalently, every eigenvalue $\lambda$ of $A$, with multiplicity $m$, corresponds to a pair of
eigenvalues $\pm\sqrt{q+\lambda}$ of $B$, each with multiplicity $m$ (or a single
eigenvalue $0$ of multiplicity $2m$ in case $\lambda=-q$). The remainder of this paper is devoted to
proving

\begin{theorem}\label{theorem1}
The graph $\Gamma=\Gamma(4,q)$ is connected except for $q\in\{2,4\}$, when the graph has
4 connected components. When $q$ is odd, the
adjacency matrix $A$ of $\Gamma$ has characteristic polynomial $\phi(x)=\det(xI_q-A)$ of the form
\[\phi(x)=(x-q(q-1))(x+q)^{(q-1)(q^2-q+1)}x^{3q(q-1)}(x-q)^{q(q-1)^2}\widetilde{\phi}(x)\]
where all roots of $\widetilde{\phi}(x)\in\Z[x]$ have the form $\lambda=-q+\varepsilon^2$ where
$|\varepsilon|\leqslant2\sqrt{q}$. Each such value $\varepsilon$ lies the ring
$\Z\bigl[2\cos\frac{2\pi}p\bigr]$, or $\Z\bigl[2\cos\frac{2\pi}9\bigr]$ if $p=3$.
\end{theorem}

A complete determination of $\phi(x)$ is given in Theorem~\ref{even} when $q$ is even,
and in Theorem~\ref{theorem2} when $q=p$ is prime. Now using Lemma~\ref{AB} we obtain

\begin{theorem}\label{B}
The graph $D(4,q)$ has eigenvalues $\pm q$, each of multiplicity 1 (unless $q\in\{2,4\}$ when each of the
eigenvalues $\pm q$ has multiplicity 4). All remaining eigenvalues have the form
$\pm\varepsilon$ where $|\varepsilon|\leqslant2\sqrt{q}$.
\end{theorem}

\noindent Once again, the eigenvalues $\varepsilon$ of Theorem~\ref{B} are cyclotomic integers
satisfying the conclusion of Theorem~\ref{theorem1}. In Theorem~\ref{theorem1} the multiplicity of the
eigenvalue~0 may actually exceed $3q(q-1)$; in particular this happens
whenever $q\equiv2\mod3$. We find explicit formulas for the actual eigenvalues, by expressing the `error'
terms $\varepsilon$ as exponential sums defined over finite fields (or over the Galois ring $GR(9,e)$ of
order $9^e=q^2$ and characteristic~9, in the case $q=3^e$). This leads to our bound
$|\varepsilon|\leqslant2\sqrt{q}$,
using the Hasse-Davenport-Weil bound when $q=p^e$, $p\geqslant5$; or the analogous bound of
Kumar, Helleseth and Calderbank~\cite{KHC} in the case $p=3$.

Our strategy for proving this result (see~\cite{B} for details) is to first realize $\Gamma$ as a Cayley graph Cay$(G,S)$ for a nonabelian
group $G$ of order $q^4$, and {\em connection set\/} $S\subset G$. (Thus $\Gamma$ has vertices labeled by elements
of $G$; and two vertices $g,g'\in G$ are adjacent in $\Gamma$ iff $g'g^{-1}\in S$). Since our graph
$\Gamma$ is undirected and connected with no loops or multiple edges, we will have
$\langle S\rangle=G$, $1\notin S$, and $g\in S$ iff $g^{-1}\!\in S$.
We then determine the number $k$ of conjugacy classes of $G$, and a complete set (up to equivalence)
of irreducible ordinary representations $\pi_i:G\to GL_{n_i}(\C)$ for $i=1,2,\ldots,k$. For each $i$, we
compute the complex $n_i\times n_i$ matrix
$\pi_i(S):=\sum_{g\in S}\pi_i(g)$.

\begin{theorem}[\cite{B,DS}; see also \cite{KS}]\label{cayley}
The characteristic polynomial of~$A$, the
adjacency matrix of~$\Gamma$, is given by
\[\phi(x)=\det(xI_{|G|}-A)=\prod_{i=1}^k\det[xI_{n_i}-\pi_i(S)]^{n_i}.\]
\end{theorem}

Note that this gives $\sum_{i=1}^k n_i^2=|G|$ eigenvalues (counting according to their respective
multiplicities) as required. In those cases where $G$ is abelian, the eigenvalues are simply the character
values $\chi_i(S)$. A similar simplification is possible when $S$ is a union of conjugacy classes of $G$, but
this does not apply in our case. When $G$ is nonabelian and the full matrices of the representations
$\pi_i$ are not explicitly known, determining the eigenvalues of $\pi_i(S)$ from the
character values alone may require substantial additional work (see~\cite{B}); but for us, the group $G$ is
sufficiently nice that explicit descriptions of
the full matrices of the representations $\pi_i$ are easily available, making our job much easier.

\section{Background on Finite Fields}\label{background}

General results on finite fields can be found in~\cite{LN}.
Let $F=\F_q$ be a field of order $q=p^e$ where $e\geqslant1$ and $p$ is prime. The absolute trace map
is
\[tr:F\to\F_p,\quad tr(a)=a+a^p+a^{p^2}+\cdots+a^{p^{e-1}}.\]
We also fix a primitive $p$-th root of unity $\zeta=\zeta_p\in\C$; here it suffices
to assume that $\zeta=e^{2\pi i/p}$. We define the {\it exponential sum\/} of an arbitrary function
$f:F\to F$ as the cyclotomic integer
\[\varepsilon_f=\sum_{a\in F}\zeta^{tr[f(a)]}\in\Z[\zeta].\]

\begin{lemma}\label{exp_sums}
For every polynomial of the form $f(t)=bt+c\in F[t]$ we have
\[\varepsilon_f=\left\{\begin{array}{ll}
0,&\hbox{if $b\neq0$;}\\
q\zeta^{tr(c)},&\hbox{otherwise.}\end{array}\right.\]
\end{lemma}

\begin{pf}
See \cite[Ch.5]{LN}.\qed
\end{pf}

\begin{lemma}\label{moments}
Let $k$ be a non-negative integer. Then
\[\sum_{a\in F}a^k=\left\{\begin{array}{ll}
-1,&\hbox{if $k=(q-1)k_1$ for some integer $k_1\geqslant1$;}\\
0,&\hbox{otherwise}.\end{array}\right.\]
\end{lemma}

\begin{pf}
See \cite[p.271]{LN}.\qed
\end{pf}

\begin{lemma}\label{count_polys}
\begin{itemize}
\item[(i)]Let $n_k$ be the number of nonzero polynomials $a_2t^2+a_1t+a_0\in F[t]$ having
exactly $k$ distinct roots in $F$. Then
\[n_0={\textstyle{\frac12}}(q-1)(q^2-q+2);\quad n_1=2q(q-1);\quad n_2={\textstyle{\frac12}}q(q-1)^2\]
and $n_k=0$ otherwise. Here $n_0+n_1+n_2=q^3-1$; and for $k=0$ we include
$\frac12 q(q-1)^2$ irreducible quadratics and $q-1$ nonzero constant polynomials.
\item[(ii)]For $q$ even, let $n_k$ be the number of nonzero polynomials $a_3t^3+a_1t+a_0\in F[t]$ having
exactly $k$ distinct\/ {\em nonzero} roots in $F$. Then
\[n_0={\textstyle{\frac13}}(q-1)(q^2+8);\quad n_1={\textstyle{\frac12}}(q-1)^2(q+4);\quad n_3={\textstyle{\frac16}}(q-1)^2(q-2)\]
and $n_k=0$ otherwise. Here $n_0+n_1+n_3=q^3-1$.
\end{itemize}
\end{lemma}

\begin{pf}
Every nonzero polynomial of degree${}\leqslant2$ with a single root has the form $a_1(t-t_1)$ or
$a_2(t-t_1)^2$, giving $n_1=2q(q-1)$. Every nonzero polynomial of degree~$2$ having two distinct roots has
the form $a_1(t-t_1)(t-t_2)$ with $a_2\neq0$ and $t_1\neq t_2$; and there are $n_2=\frac12 q(q-1)^2$
such polynomials. This leaves $n_0=q^3-1-n_1-n_2=\frac12(q-1)(q^2-q+2)$, and the remaining assertions
of (i) follow.

Now suppose $q$ is even, and consider a nonzero polynomial $f(t)=a_3t^3+a_1t+a_0\in F[t]$.
If $f(t)=a_3(t+t_1)(t+t_2)(t+t_3)$ then $t_1+t_2+t_3=0$; so in characteristic~2, the number of distinct
nonzero roots must
be 0, 1 or 3. There are $n_3=\frac16(q-1)^2(q-2)$ nonzero polynomials of the form
$f(t)=a_3(t+t_1)(t+t_2)(t+t_1+t_2)$ where $t_1,t_2$ are nonzero and distinct. There are $(q-1)^2$ cubics
of the form $a_3t(t+t_1)^2$
where $t_1\neq0$; and by (i), there are $\frac12q(q-1)^2$ cubics of the form
$(t+t_1)\bigl(a_3t^2+a_3 t_1 t+\frac{a_0}{t_1}\bigr)$ for which $t_1\neq0$ and the quadratic factor is
irreducible. These, together with the $(q-1)^2$ polynomials $a_1t+a_0$ having $a_1,a_0\neq0$, give
\[n_1=(q-1)^2+{\textstyle{\frac12}}q(q-1)^2+(q-1)^2={\textstyle{\frac12}}(q-1)(q^2+4).\]
This leaves
\[n_0=q^3-1-n_1-n_3={\textstyle{\frac12}}(q-1)(q^2+8).\]
One checks that this includes $\frac13(q-1)(q^2+2)$ irreducible cubics of the required form, together
with $q-1$ polynomials of the form $a_1 t$ with $a_1\neq0$, and $q-1$ nonzero constant polynomials.\qed
\end{pf}

\section{A Regular Group of Automorphisms of $\Gamma$}\label{main_proof}

For all $t,u,v,w\in F$ we define the matrix
\[g=g(t,u,v,w)=\left[\begin{array}{ccccc}1&t&u&v{+}tu&w\\
&1&0&-u&0\\
&&1&t&0\\
&&&1&0\\
&&&&1\end{array}\right].\]
These $q^4$ matrices form a subgroup $G<GL_5(F)$ acting regularly on points via
\[(1,p_1,p_2,p_3,p_4)\mapsto(1,p_1,p_2,p_3,p_4)g(t,u,v,w).\]
which can be written simply as
\[P\mapsto Pg\]
after a slight abuse of notation by which we identify
\[P=P(p_1,p_2,p_3,p_4)=(1,p_1,p_2,p_3,p_4).\]
One checks that this action preserves collinearity of points, and so gives a group of automorphisms
of $\Gamma$ which is regular on the vertices. Thus $\Gamma$ is a Cayley graph Cay$(G,S)$
for the set of $q(q-1)$ elements
\begin{align*}
S&{}=\{g\in G\,:\,P(0,0,0,0)g\hbox{\ is (distinct from and) collinear with\ }P(0,0,0,0)\}\\
&{}=\{g(t,rt,-rt^2,r^2 t)\,:\,r,t\in F,\;t\neq0\}.\end{align*}
The commutator of two typical elements of $G$ is
\[[g(t,u,v,w),g(t',u',v',w')]=g(0,0,2t'u{-}2tu',0).\]
At this point we must consider separately the cases $q$ even and $q$ odd, for which $G$ is abelian or
nonabelian, respectively.

\section{The case $q$ even}\label{q_even}

In this section we suppose $q$ is even, so that
\[g(t,u,v,w)g(t',u',v',w')=g(t+t',u+u',v+v',w+w').\]
In this case $G$ is elementary abelian, with $q^4$ irreducible linear characters
\[\chi_{\alpha,\beta,\gamma,\eta}\bigl(g(t,u,v,w)\bigr)=(-1)^{tr(\alpha t+\beta u+\gamma v+\eta w)},
\quad\alpha,\beta,\gamma,\eta\in F.\]

\begin{theorem}\label{even}
Suppose $q$ is even. Then the characteristic polynomial of $A$, the incidence matrix of $\Gamma$, is
\begin{align*}
\phi(x)&{}=\det(xI_{q^4}-A)\\
&{}=\bigl(x-q(q-1)\bigr)(x-3q)^{q(q-1)^2(q-2)/24}(x-q)^{q(q-1)^2(q+4)/4}\\
&\quad{}\times x^{(q-1)(q^3+8q+3)/3}(x+q)^{3q(q-1)^2(q+2)/8}.
\end{align*}
The graph $\Gamma$ is connected for $q\geqslant8$; while for $q\in\{2,4\}$, $\Gamma$ has 4 connected
components.
\end{theorem}

\begin{pf}By Theorem~\ref{cayley} and Lemma~\ref{exp_sums}, we have
\begin{align*}
\phi(x)&{}=\prod_{\alpha,\beta,\gamma,\eta\in F}\Bigl(x-\sum_{\myatop{r,t\in F}{t\neq0}}(-1)^{tr(\alpha t+\beta rt+\gamma rt^2+\eta r^2 t)}\Bigr)\\
&{}=(x-q(q-1))(x+q)^{q-1}\mskip-12mu\prod_{\myatop{\alpha,\beta,\gamma,\eta\in F}{(\beta,\gamma,\eta)\neq(0,0,0)}}\mskip-12mu\Bigl(x-\sum_{\myatop{r,t\in F}{t\neq0}}(-1)^{tr(\alpha t+\beta rt+\gamma rt^2+\eta r^2 t)}\Bigr).
\end{align*}
Now using the fact that the map $F\to F$, $r\mapsto r^2$ is an automorphism (in particular bijective and
trace-preserving),
\begin{align*}
\sum_{\myatop{r,t\in F}{t\neq0}}(-1)^{tr(\alpha t+\beta rt+\gamma rt^2+\eta r^2 t)}
&{}=\sum_{0\neq t\in F}(-1)^{tr(\alpha t)}\sum_{r\in F}(-1)^{tr[(\beta^2 t+\gamma^2 t^3+\eta)r^2 t]}\\
&{}=q\mskip-10mu\sum_{\myatop{0\neq t\in F}{\beta^2t+\gamma^2t^3=\eta}}\mskip-10mu(-1)^{tr(\alpha t)}.
\end{align*}
After re-indexing via $(\beta,\gamma,\eta)\mapsto(\beta^{1/2},\gamma^{1/2},\eta)$,
\begin{align*}
\phi(x)&{}=(x-q(q-1))(x+q)^{q-1}\mskip-12mu\prod_{(\beta,\gamma,\eta)\neq(0,0,0)}\prod_{\alpha}
\,\Bigl(x-q\mskip-10mu\sum_{\myatop{0\neq t\in F}{\beta t+\gamma t^3=\eta}}\mskip-10mu(-1)^{tr(\alpha t)}\Bigr).
\end{align*}
If the polynomial $f(t)=\gamma t^3+\beta t+\eta\in F[t]$ has a unique nonzero root $t_1\in F$,
then the map $F\to\F_2$, $\alpha\mapsto tr(\alpha t_1)$ takes each of the values in $\{0,1\}$
exactly $\frac{q}2$ times, in which case
\[\prod_\alpha\,\Bigl(x-q\mskip-10mu\sum_{\myatop{0\neq t\in F}{\beta t+\gamma t^3=\eta}}\mskip-10mu(-1)^{tr(\alpha t)}\Bigr)=(x^2-q^2)^{q/2}.\]
Similarly, if $f(t)$ (as above) has three distinct nonzero
roots $t_1,t_2,t_3\in F$, then $t_1+t_2+t_3=0$ and the map $F\to\F_2^3$,
$\alpha\mapsto(tr(\alpha t_1),tr(\alpha t_2),tr(\alpha t_3))$ attains each of the triples
$(0,0,0),(1,1,0),(1,0,1),(0,1,1)$ exactly $\frac{q}4$ times, in which case
\[\prod_\alpha\,\Bigl(x-q\mskip-10mu\sum_{\myatop{0\neq t\in F}{\beta t+\gamma t^3=\eta}}\mskip-10mu(-1)^{tr(\alpha t)}\Bigr)=(x-3q)^{q/4}(x+q)^{3q/4}.\]
Thus
\[\phi(x)=(x-q(q-1))(x+q)^{q-1}x^{n_0}\bigl[(x^2-q^2)^{q/2}\bigr]^{n_1}
\bigl[(x-3q)^{q/4}(x+q)^{3q/4}\bigr]^{n_3}\]
where $n_k$ is given by Lemma~\ref{count_polys}(ii). Simplification yields the formula claimed for $\phi(x)$.
Now we simply read off the multiplicity of the largest eigenvalue to obtain the number of connected
components of $\Gamma$ (see e.g.~\cite[Prop.1.3.8]{BH}).\qed
\end{pf}

\section{The case $q$ is odd}\label{q_odd}

Here and {\em for the remainder of this paper, we take $q$ to be odd.\/}
From the general formula for commutators in $G$ given at the end of Section~\ref{main_proof}, we
deduce the commutator subgroup and centre
\[G'=\{g(0,0,u,0):u\in F\},\quad Z=Z(G)=\{P(0,0,v,w):v,w\in F\};\]
also the centralizer of a noncentral element (i.e.\ with $(t,u)\neq(0,0)$) is a subgroup
\[C_G\bigl(g(t,u,v,w)\bigr)=\bigl\{g(ct,cu,v',w'):c,v',w'\in F\bigr\}\]
of order $q^3$. So $G$ has $q^3+q^2-q$ conjugacy classes ($q^2$ of size 1, and $q^3-q$ of size
$q$). There are $|G/G'|=q^3$ linear characters of~$G$, given by
\[\chi_{\alpha,\beta,\gamma}\bigl(g(t,u,v,w)\bigr)=\zeta^{tr(\alpha t+\beta u+\gamma w)}\]
where $\alpha,\beta,\gamma\in F$. As in Section~\ref{background}, $\zeta=\zeta_p$ is a complex $p$-th root of
unity and $tr:F\to\F_p$ is the trace map.
The remaining irreducible characters of $G$ may be found by inducing linear
characters of a subgroup of order $q^3$ (thus yielding monomial representations of degree $q$); but
guided by a little hindsight, we will instead directly exhibit the
missing representations and show that they are irreducible and distinct. For each pair $\alpha,\beta\in F$ with
$\alpha\neq0$, we define $M_{\alpha,\beta}:G\to GL_q(\C)$ by
\[M_{\alpha,\beta}(g(t,u,v,w))=\bigl[\zeta^{tr[\alpha(v-2iu)+\beta w]}\delta_{i+t,j}\bigr]_{i,j\in F}\]
using the Kronecker delta notation $\delta_{i,j}=0$ or $1$ according as $i,j\in F$ either differ or coincide.
It is routine to check that $M_{\alpha,\beta}(g)M_{\alpha,\beta}(g')=M_{\alpha,\beta}(gg')$ for all $g,g'\in G$,
and $M_{\alpha,\beta}(g(0,0,0,0))=I_q$; so $M_{\alpha,\beta}$ is a
representation of degree~$q$. The associated character is found to be
\[\psi_{\alpha,\beta}\bigl(g(t,u,v,w)\bigr)=tr\,M_{\alpha,\beta}\bigl(g(t,u,v,w)\bigr)
=\left\{\begin{array}{ll}
\zeta^{tr(\alpha v+\beta w)}q,&\hbox{if $t=u=0$;}\\
0,&\hbox{otherwise}
\end{array}\right.\]
using the fact that $\alpha\neq0$. These $q^2-q$ characters of $G$ are irreducible and inequivalent since
\begin{align*}[\psi_{\alpha,\beta},\psi_{\alpha',\beta'}]_G
&{}=\frac1{q^4}\sum_{t,u,v,w\in F}\psi_{\alpha,\beta}\bigl(g(t,u,v,w\bigr)\overline{\psi_{\alpha',\beta'}\bigl(g(t,u,v,w)\bigr)}\\
&{}=\frac1{q^4}\sum_{v,w\in F}\zeta^{tr[(\alpha-\alpha')v+(\beta-\beta')w]}q^2\\
&{}=\left\{\begin{array}{ll}
1,&\hbox{if $(\alpha,\beta)=(\alpha',\beta')$;}\\
0,&\hbox{otherwise.}
\end{array}\right.\end{align*}
These are also distinct from the characters $\chi_{\alpha,\beta,\gamma}$ and so we have the complete list of
$q^3+q^2-q$ irreducible characters of $G$.

Now by Theorem~\ref{cayley}, the adjacency matrix $A$ of $\Gamma$ has characteristic polynomial
\[\phi(x)=\det(xI_{q^4}-A)=\prod_{\alpha,\beta,\gamma\in F}(x-\chi_{\alpha,\beta,\gamma}(S))
\prod_{\myatop{\alpha,\beta\in F}{\alpha\neq0}}\det[xI_q-M_{\alpha,\beta}(S)]^q.\]
Those eigenvalues of $A$ obtained from the linear characters of $G$ are
\[\chi_{\alpha,\beta,\gamma}(S)=\sum_{g\in S}\chi_{\alpha,\beta,\gamma}(g)
=\sum_{\myatop{r,t\in F}{t\neq0}}\zeta^{tr[(\alpha+\beta r+\gamma r^2)t]}=(m_{\alpha,\beta,\gamma}-1)q\]
where $m_{\alpha,\beta,\gamma}$
is the number of values $r\in F$ such that $\alpha+\beta r+\gamma r^2=0$. By Lemma~\ref{count_polys}(ii), the first $q^3$ factors of $\phi(x)$ are
\begin{align*}
(x-q(q-1))\mskip-12mu&\prod_{(\alpha,\beta,\gamma)\neq(0,0,0)}\mskip-4mu\Bigl(x-(m_{\alpha,\beta,\gamma}-1)q\Bigr)\\
&{}=(x-q(q-1))(x-q)^{n_2}x^{n_1}(x+q)^{n_0}\\
&{}=(x-q(q-1))(x-q)^{q(q-1)^2/2}x^{2q(q-1}(x+q)^{(q-1)(q^2-q+2)/2};
\end{align*}
thus the characteristic polynomial $\phi(x)=\det(xI_{q^4}-A)$ has the form
\begin{align*}
\phi(x)={}&x^{2q(q-1)}(x-q(q-1))(x-q)^{q(q-1)^2/2}(x+q)^{(q-1)(q^2-q+2)/2}\\
&{}\times\prod_{\myatop{\alpha,\beta\in F}
{\alpha\neq0}}\mskip-5mu\bigl[\det(xI_q-M_{\alpha,\beta}(S))\bigr]^q.\end{align*}
Now for $\alpha\neq0$,
\begin{align*}
M_{\alpha,\beta}(S)&{}=\sum_{g\in S}M_{\alpha,\beta}(g)
=\sum_{\myatop{r,t\in F}{t\neq0}}\bigl[\zeta^{tr[\beta r^2 t-\alpha(t+2i)rt]}\delta_{i+t,j}\bigr]_{i,j\in F}\\
&{}=\Bigl[\,{\textstyle\sum\limits_{r\in F}}\zeta^{tr[\beta r^2(j-i)-\alpha r(j^2-i^2)]}\Bigr]_{i,j\in F}-qI_q\\
&{}=U_{\alpha,\beta}U_{\alpha,\beta}^*-qI_q\end{align*}
where `$*$' denotes conjugate-transpose, and we have introduced the $q\times q$ complex matrices
\[U_{\alpha,\beta}=\bigl[\zeta^{tr(\alpha i^2 j-\beta ij^2)}\bigr]_{i,j\in F}\,.\]
We first treat the cases $\beta=0\neq\alpha$ for which we obtain
\[M_{\alpha,0}(S)=\Bigl[\,{\textstyle\sum\limits_{r\in F}}\zeta^{tr[\alpha r(i^2-j^2)]}\Bigr]_{i,j\in F}-qI_q.\]
\vskip-3pt\noindent
Denoting by $\{e_r\}_{r\in F}$ the standard basis of $\C^F=\C^q$, we find a new basis consisting of
eigenvectors of $M_{\alpha,0}(S)$ as follows:
\begin{itemize}
\item$\frac12(q-1)$ eigenvectors of the form $e_r+e_{-r}$ where $0\neq r\in F$, each with
eigenvalue $q$;
\item$\frac12(q-1)$ eigenvectors of the form $e_r-e_{-r}$ as $r$ ranges over a set of
representatives of the distinct nonzero pairs $\{r,-r\}$ in $F$. Each such vector has
eigenvalue $-q$;
\item$M_{\alpha,0}e_0=0$.
\end{itemize}

\smallskip\noindent
After including the factors
\[\prod_{0\neq\alpha\in F}\bigl[\det(xI_q-M_{\alpha,0}(S))\bigr]^q
=x^{q(q-1)}(x^2-q^2)^{q(q-1)^2/2},\]
we update our formula for the characteristic polynomial of $A$ as
\begin{align*}
\phi(x)={}&x^{3q(q-1)}(x-q(q-1))(x-q)^{q(q-1)^2}(x+q)^{(q-1)(2q^2-2q+1)}\\
&{}\times\prod_{\myatop{\alpha,\beta\in F}
{\alpha\beta\neq0}}\mskip-5mu\bigl[\det(xI_q-M_{\alpha,\beta}(S))\bigr]^q.\end{align*}
Finally we describe the remaining $q^2(q-1)^2$ eigenvalues of $A$ arising from
$M_{\alpha,\beta}(S)$ for $\alpha\beta\neq0$.

\begin{lemma}\label{lemma1}
For any nonzero elements $c,d\in F$ the matrix
$M_{\alpha,\beta}(S)$ is similar to $M_{c^2 d\alpha,cd^2\beta}(S).$
\end{lemma}

\begin{pf}Re-indexing rows and columns of $U_{\alpha,\beta}$ via $(i,j)\mapsto(ci,dj)$, we see that
\[U_{c^2 d\alpha,cd^2\beta}=P_c U_{\alpha,\beta}P_d^T\]
where $P_c$ and $P_d$ are $q\times q$ permutation matrices,
and so $M_{c^2 d\alpha,cd^2\beta}=P_c M_{\alpha,\beta}P_c^T$.\break\qed
\end{pf}

\begin{corollary}\label{cor1}
If $q\not\equiv1\mod3$ then for all nonzero $\alpha,\beta\in F$, $M_{\alpha,\beta}(S)$ is
unitarily similar to $M_{1,1}(S)$. If $q\equiv1\mod3$ then there are at most three similarity classes of matrices
$M_{\alpha,\beta}(S)$ with $\alpha\beta\neq0$, represented by $M_{1,1}(S)$, $M_{1,\omega}(S)$ and
$M_{1,\omega^2}(S)$ where $\omega\in F$ is a primitive root.
\end{corollary}

\begin{pf}If $q\not\equiv1\mod3$, then every element of $F$ has a cube root in~$F$; so let $c\in F$ be any cube root
of $\beta/\alpha^2$ and take $d=\alpha c/\beta$. Then $M_{\alpha,\beta}(S)$ is similar to
$M_{1,1}(S)$ by Lemma~\ref{lemma1}. The second conclusion follows similarly.\qed
\end{pf}

If $q\equiv2\mod3$ and $\alpha\beta\neq0$, then $M_{\alpha,\beta}(S)$ is similar to $M_{3,3}$,
since by Corollary~\ref{cor1}, both matrices are similar to $M_{1,1}$.
In this case $U_{3,3}=\bigl[\zeta^{tr(3i^2 j-3ij^2)}\bigr]_{i,j\in F}$ is unitarily similar to
\[\tilde{U}:=D^*U_{3,3}D=\bigl[\zeta^{tr[(j-i)^3]}\bigr]_{i,j\in F}\]
where $D$ is a diagonal matrix with diagonal entries $\zeta^{tr(i^3)}$ for $i\in F$.
In this case the $q$ vectors $v_c=\bigl(\zeta^{tr(ci)}\bigr)_{i\in F}$ for $c\in F$ form a basis of $\C^F$
consisting of eigenvectors of $\tilde{U}$; indeed
\[(\tilde{U}v_c)_i=\sum_{j\in F}\zeta^{tr[(j-i)^3]}\zeta^{tr(cj)}
=\sum_{j\in F}\zeta^{tr(j^3)}\zeta^{tr(cj+ci)}=\Bigl(\sum_{j\in F}\zeta^{tr(j^3+cj)}\Bigr)\zeta^{tr(ci)}\]
so that $\tilde{U}v_c=\varepsilon_f v_c$ where
\[\varepsilon_f=\sum_{r\in F}\zeta^{tr\,f(r)};\quad f(x)=x^3+cx\in F[x].\]
Note that $\varepsilon_f\in\Z[\zeta]$ satisfies $\overline{\varepsilon_f}=\varepsilon_f$ since $f(-x)=-f(x)$;
thus $\varepsilon_f\in\Z[\zeta+\overline{\zeta}]=\Z\bigl[2\cos\frac{2\pi}p\bigr]$
(see \cite[Prop.2.16]{Wa}). Also $M_{3,3}=D^*\tilde{U}\tilde{U}^*D-qI_q$ has eigenvalues
$\varepsilon_f^2-q$. The Weil bound (see e.g.~\cite{W}, \cite[p.223]{LN}) gives
$|\varepsilon_f|\leqslant2\sqrt{q}$ as required.
The all-ones eigenvector $v_0$ has eigenvalue $\sum_{j\in F}\zeta^{tr(j^3)}=0$ since
$f(x)=x^3$ is a permutation of $F$; so for $q\equiv2\mod3$ we obtain
\begin{align*}
\phi(x)={}&x^{3q(q-1)}(x-q(q-1))(x-q)^{q(q-1)^2}(x+q)^{(q-1)(2q^2-2q+1)}\\
&{}\times\prod_{0\neq c\in F}\mskip-4mu(x+q-\varepsilon_{t^3+ct}^2)^{q(q-1)}.\end{align*}

When $q\equiv1\mod3$ we work just a little harder. Let $a=\frac1{3\alpha\beta}$ and use the identity
\[a(\beta j-\alpha i)^3=a\beta^3 j^3-\beta ij^2+\alpha i^2 j-a\alpha^3 i^3\]
to see that $U_{\alpha,\beta}=N\tilde{U}N'$ where
\[\tilde{U}=\bigl[\zeta^{tr[a(j-i)^3]}\bigr]_{i,j\in F}\]
and the unitary matrices $N$ and $N'$ are given by
\[N=\bigl[\zeta^{tr(aj^3)}\delta_{\alpha i,j}\bigr]_{i,j\in F}\,;\quad
N'=\bigl[\zeta^{-tr(ai^3)}\delta_{i,\beta j}\bigr]_{i,j\in F}\,.\]
Now the vectors $v_c$ (as above) are eigenvectors of $\tilde{U}$ since
\[(\tilde{U}v_c)_i=\sum_{j\in F}\zeta^{tr[a(j-i)^3]}\zeta^{tr(cj)}
=\sum_{j\in F}\zeta^{tr(aj^3)}\zeta^{tr(cj+ci)}=\Bigl(\sum_{j\in F}\zeta^{tr(aj^3+cj)}\Bigr)\zeta^{tr(ci)}\]
\vskip-5pt\noindent
with corresponding eigenvalue $\varepsilon_f$ where $f(t)=at^3+ct\in F[t]$. Now
$M_{\alpha,\beta}=N\tilde{U}\tilde{U}^*N^*-qI_q$ has eigenvalues $\varepsilon_f^2-q$, and
\[\prod_{\myatop{\alpha,\beta\in F}{\alpha\beta\neq0}}\det(xI_q-M_{\alpha,\beta}(S))
=\prod_{\myatop{a,c\in F}{a\neq0}}(x+q-\varepsilon_{at^3+t}^2)^{q-1}\]
since as $\alpha$ and $\beta$ range over the nonzero elements of $F$, the coefficient
$a=\frac1{3\alpha\beta}$ falls in each of the three multiplicative cosets of the cubes equally often,
and we recall Corollary~\ref{cor1}. So for $q\equiv1\mod3$ we obtain
\begin{align*}
\phi(x)={}&x^{3q(q-1)}(x-q(q-1))(x-q)^{q(q-1)^2}(x+q)^{(q-1)(q^2-q+1)}\\
&{}\times\prod_{\myatop{a,c\in F}{a\neq0}}(x+q-\varepsilon_{at^3+ct}^2)^{q(q-1)}.\end{align*}
As before, $\varepsilon_f\in\Z\bigl[2\cos\frac{2\pi}p\bigr]$ and $|\varepsilon_f|\leqslant2\sqrt{q}$.

Finally, suppose $q=3^e$ so that $p=3$. By the Lemma, $U_{\alpha,\beta}$ is unitarily similar to $U_{1,1}$.
Unlike the cases $q\not\equiv0\mod3$, in this case the eigenvalues of $U_{1,1}$ do not lie in $\Q[\zeta]$;
rather they lie in $\Q[\xi]$ where we abbreivate $\xi=\zeta_9=e^{2\pi i/9}$, chosen so that
$\xi^3=\zeta=\zeta_3$. This can be seen even in the case $q=3$ where
\[U_{1,1}=\left[\begin{array}{ccc}1&1&1\\ 1&1&\zeta\\ 1&\zeta^2&1\end{array}\right]\]
whose eigenvalues are
\[1+\xi^4+\xi^5,\quad1+\xi^2+\xi^7=1-\xi+\xi^2-\xi^4,\quad
1+\xi+\xi^8=1+\xi-\xi^2-\xi^5.\]
Likewise, the eigenvalues of $M_{1,1}$ lie in $\Q[\xi]$ but not in $\Q[\zeta]$. However we see that the
eigenvalues of $U_{1,1}$ are expressible as exponential sums defined over Galois rings;
see e.g.~\cite{KHC,M}. Let $R=GR(9,e)$ be the Galois ring of order
$9^e=q^2$ and characteristic~9. The ring $R$ enjoys the following properties:
\smallskip
\begin{itemize}
\item{}$R$ is a commutative ring with a maximal ideal $3R$ consisting of all zero divisors in $R$,
and the quotient ring is $R/3R\cong F=\F_q$.
\item{}The units of $R$ form a multiplicative group $R^\times$ consisting of all elements not in $3R$.
This group has a multiplicative subgroup $\{1,\beta,\beta^2,\ldots,\beta^{q-2}\}$ of order $q-1$.
\item{}Every element $x\in R$ has a unique 3-adic expansion $x=x_0+3x_1$ where $x_0,x_1\in\T$.
where we define $\T=\{0,1,\beta,\beta^2,\ldots,\beta^{q-2}\}$. In particular, $\T$ is a set of representatives
of the cosets $R/3R\cong\F_q$.
\item{}The trace map $tr:R\to\Z/9\Z$ is defined by
\[tr(x_0+3x_1)=(x_0+x_0^3+x_0^9+\cdots+x_0^{3^{e-1}})+3(x_1+x_1^3
+x_1^9+\cdots+x_1^{3^{e-1}})\]
where $x_0,x_1\in\T$. After reducing both domain and range modulo~$3R$, this gives the usual
absolute trace map $\F_q\to\F_3$.
\end{itemize}
\smallskip

After replacing $F=\F_q$ by $\T$ as index set for entries of our vectors and matrices, we may rewrite our
basis of $\C^q$ as $v_c=\bigl(\xi^{3tr(ci)}\bigr)_{i\in\T}$ for $c\in\T$, and
\[U_{1,1}=\bigl[\xi^{3tr(i^2 j-ij^2)}\bigr]_{i,j\in\T}\]
which is unitarily similar to
\[\tilde{U}=\bigl[\xi^{tr[(j-i)^3]}\bigr]_{i,j\in\T}\]
after conjugating by the unitary diagonal matrix $D=\bigl[\xi^{tr(i^3)}\delta_{i,j}\bigr]_{i,j\in\T}\,$.
Now
\[(\tilde{U}v_c)_i=\sum_{j\in\T}\xi^{tr[(j-i)^3]}\xi^{3tr(cj)}
=\sum_{j\in\T}\xi^{tr(j^3)}\xi^{3tr(cj+ci)}=\Bigl(\sum_{j\in\T}\xi^{tr(j^3+3cj)}\Bigr)\xi^{3tr(ci)}\]
so that $\tilde{U}v_c=\varepsilon_f v_c$ where
\[\varepsilon_f=\sum_{i\in\T}\xi^{tr\,f(i)};\quad f(t)=t^3+3ct\in R[t],\;c\in\T.\]
The `weighted degree' of $f(t)$, as defined in~\cite[p.458]{KHC}, is $d=3$; and as shown in~\cite{KHC},
the Weil bound $|\varepsilon_f|\leqslant(d-1)\sqrt{q}=2\sqrt{q}$ holds. Once again by
\cite[Prop.2.16]{Wa} we have $\varepsilon_f\in\Z\bigl[2\cos\frac{2\pi}9\bigr]$.
With the notation above, we have
\begin{align*}
\phi(x)={}&x^{3q(q-1)}(x-q(q-1))(x-q)^{q(q-1)^2}(x+q)^{(q-1)(q^2-q+1)}\\
&{}\times\prod_{c\in\T}(x+q-\varepsilon_{t^3+3ct}^2)^{q(q-1)^2}.\end{align*}

\def\inv{{}^{\scriptscriptstyle-1}{\mskip-2mu}}
\section{Exact Spectra over Prime Fields}\label{exact}

It is possible to refine Theorem~\ref{theorem1} to express $\phi(x)$ precisely. Here, however, we state such a result (Theorem~\ref{theorem2} below) only for the case
$q=p$ is prime; in the general case $q=p^e$ with $e\geqslant2$, counting multiplicities is more
technical and Theorem~\ref{theorem1} is probably adequate for any intended applications. In the following
we denote $f\inv(r)=\{a\in F:f(a)=r\}$.

\begin{lemma}\label{fibres}
Let $f,g:F\to F$ be two functions where the field $F=\F_p$ has prime order~$p$. Then
$\varepsilon_f=\varepsilon_g$ iff $|f\inv(r)|=|g\inv(r)|$ for all $r\in F$.
\end{lemma}

\begin{pf}
If $\varepsilon_f=\varepsilon_g$ then
\[0=\varepsilon_f-\varepsilon_g=\sum_{r\in F}\bigl(|f\inv(r)|-|g\inv(r)|\bigr)\zeta^r.\]
Since the minimal polynomial of $\zeta$ over $\Q$ is the cyclotomic
polynomial $1+x+x^2+\cdots+x^{p-1}\in\Z[x]$, there exists $n\in\Z$ such that
$|f\inv(r)|-|g\inv(r)|=n$ for all $r\in F$. Since $\sum_{r\in F}|f\inv(r)|=\sum_{r\in F}|g\inv(r)|=p$,
we must have $n=0$. The converse is clear.\qed
\end{pf}

\begin{lemma}\label{lambda}
Let $f(t)\in F[t]$ and $\tilde{f}(t)=f(\lambda t)\in F[t]$ where $0\neq\lambda\in F$. Then
$\varepsilon_{\tilde{f}}=\varepsilon_f$.
\end{lemma}

\begin{pf}
Straightforward.\qed
\end{pf}

Let $C$ be the set of $p(p-1)$ cubic polynomials of the form $f(t)=at^3+ct\in F[t]$. By
Lemma~\ref{lambda}, the corresponding exponential sums $\varepsilon_f$ for $f\in C$ are not all distinct.
We next find a set of representatives $\tilde{C}\subset C$ giving rise to distinct exponential sums;
that is, for each $f\in C$ there is a unique $g\in\tilde{C}$ such that $\varepsilon_f=\varepsilon_g$.

\begin{lemma}\label{cubic_sums}
Let $f(t)=at^3+ct\in F[t]$ where the field $F=\F_p$ has prime order $p\geqslant5$, and
$a\neq0$.
\begin{itemize}
\item[(i)]For $p\equiv2\mod3$, we have $\varepsilon_f=0$ iff $c=0$. We may take
$\tilde{C}=\{t^3+\tilde{c}t:\tilde{c}\in F\}$; and $\varepsilon_f=\varepsilon_{t^3+\tilde{c}t}$
iff $\tilde{c}=a^{-1/3}c$. Here $|\tilde{C}|=p$.
\item[(ii)]For $p\equiv1\mod3$, we have $\varepsilon_f\neq0$. Let $\omega\in F$ be a
primitive element, i.e.\ a generator of the multiplicative group $F^\times$.
We may take $\tilde{C}=\{t^3,\omega t^3,\omega^2 t^3\}\cup\{\tilde{a}t^3+t:0\neq\tilde{a}\in F\}$.
Here $|\tilde{C}|=p+2$. If $c\neq0$ then $\varepsilon_f=\varepsilon_{\tilde{a}t^3+1}$ where $\tilde{a}=a/c^3$.
If $c=0$ then $\varepsilon_f=\varepsilon_{\omega^i t^3}$ where $i\in\{0,1,2\}$ is uniquely
determined by $(a/\omega^i)^{(p-1)/3}=1$.
\end{itemize}
\end{lemma}

\begin{pf}(i) First suppose $p\equiv2\mod3$, so that $p=3m-1$ where $m\geqslant2$. In this case
every element $a\in F$ has a unique cube root $a^{1/3}=a^{2m-1}$.
By Lemma~\ref{moments},
\[\sum_{s\in F}|f\inv(s)|s^m=\sum_{r\in F}(ar^3+cr)^m=\sum_{i=0}^m{m\choose i}a^{m-i}c^i\sum_{r\in F}r^{3m-2i}=-ma^{m-1}c.\]
If $\varepsilon_f=\varepsilon_g$ where $g(t)=t^3+\tilde{c}t$ then $|f\inv(s)|=|g\inv(s)|$ for all $s\in F$
by Lemma~\ref{fibres}, so $-ma^{m-1}c=-m\tilde{c}$ and $\tilde{c}=a^{m-1}c=a^{-1/3}c$;
but conversely, if $\tilde{c}=a^{-1/3}c$ then $f(t)=g(a^{1/3}t)$ so by
Lemma~\ref{lambda}, $\varepsilon_f=\varepsilon_g$.

If $c=0$ then $f(t)=at^3$ defines a permutation of $F$, so $\varepsilon_f=0$.
\smallskip

(ii) Now suppose $p\equiv1\mod3$, and write $p=3m+1$ where $m\geqslant2$.
Again by Lemma~\ref{moments},
\[\sum_{s\in F}|f\inv(s)|s^{m}=\sum_{r\in F}(ar^3+cr)^m
=\sum_{i=0}^m{m\choose i}a^{m-i}c^i\sum_{r\in F}r^{3m-2i}=-a^m\]
and
\begin{align*}
\sum_{s\in F}|f\inv(s)|s^{m+2}&=\sum_{r\in F}(ar^3+cr)^{m+2}
=\sum_{i=0}^{m+2}{m+2\choose i}a^{m+2-i}c^i\sum_{r\in F}r^{3m+6-2i}\\
&=-{m+2\choose3}a^{m-1}c^3.
\end{align*}
If $\varepsilon_f=\varepsilon_g$ where $g(t)=\tilde{a}t^3+\tilde{c}t$, then as in (i), it follows that
$\tilde{a}^m=a^m$, $\tilde{a}^{m-1}\tilde{c}^3=a^{m-1}c^3$ and $\tilde{a}c^3=a\tilde{c}^3$.
We consider two cases:
\smallskip
\begin{itemize}
\item Suppose $c=0$; then $\tilde{c}=0$. We may write $a=\omega^{3d+i}$ where $d\in\Z$ and
$i\in\{0,1,2\}$, so that $(\omega^i)^m=a^m$.
Now if $g\in\tilde{C}$ satisfies $\varepsilon_g=\varepsilon_f$, we must have $g(t)=\omega^i t^3$.
Conversely, $\omega^i t^3=f(\omega^{-d}t)$; so $g(t)=\omega^i t^3$ is the unique $g\in\tilde{C}$
satisfying $\varepsilon_g=\varepsilon_f$.
\item Suppose $c\neq0$; then $\tilde{c}\neq0$. If $g\in\tilde{C}$ satisfies $\varepsilon_g=\varepsilon_f$,
we must have $g(t)=\tilde{a}t^3+t$ where $\tilde{a}=c^{-3}a$. Conversely, $c^{-3}at^3+t=f(c^{-1}t)$;
so $g(t)=c^{-3}at^3+t$ is the unique $g\in\tilde{C}$ satisfying $\varepsilon_g=\varepsilon_f$.\qed
\end{itemize}
\end{pf}

In order to determine the exact spectrum of $\Gamma(4,p)$, we need to know not only when the values
$\varepsilon_f$ are distinct, but actually when the values $\varepsilon_f^2$ are distinct.
As preparation, we need the following.

\begin{lemma}\label{triple}
Let $F=\F_p$ where $p\geqslant5$ is prime, and suppose there exists a polynomial of the form
$f(x)=x^3+cx\in F[t]$ having $|f\inv(s)|\leqslant2$ for all $s\in F$. Then one of the following holds:
\begin{itemize}
\item[(i)]$p\equiv2\mod3$, and $c=0$;
\item[(ii)]$p=5$ and $c\in\{2,3\}$; or
\item[(iii)]$p=7$ and $c\in\{1,2,4\}$.
\end{itemize}
\end{lemma}

\begin{pf}We will assume $p\geqslant11$ since the cases $p\in\{5,7\}$ may be easily checked by explicit
computation.
First observe that $|f\inv(s)|=2$ for at most two values of $s\in F$; this is because for any such value of
$s$, $f(t)-s=t^3+ct-s$ and $f'(t)=3t^2+c$ have a linear factor in common, forcing $s^2=-4c^3/27$.

For each $k\in\{0,1,2,3\}$, let $n_k=|\{s\in F\,:\,|f\inv(s)|=k\}|$. By hypothesis, $n_3=0$; and we have just shown that
$n_2\leqslant2$. Elementary counting arguments give $n_0+n_1+n_2=p=n_1+2n_2$. Exactly three possibilities must be considered.

Case (i): $(n_0,n_1,n_2)=(0,p,0)$. In this case $f$ is a permutation polynomial; but then $\varepsilon_f=0$ and Lemma~\ref{cubic_sums} gives $p\equiv2\mod3$, and $f(t)=t^3$.

Case (ii): $(n_0,n_1,n_2)=(1,p{-}2,1)$. Here $|f\inv(a)|=2$ and $|f\inv(b)|=0$ for some $a,b\in F$, and
$|f\inv(s)|=1$ for all other values $s\in F$. Since $f$ has degree $3<p-1$, $0=\sum_{t\in F}f(t)=a-b$,
contradicting Lemma~\ref{moments}.

Case (iii): $(n_0,n_1,n_2)=(2,p{-}4,2)$. Here there exist distinct values $a_0,a_1$, $a_2,a_3\in F$ such that
\[|f\inv(s)|=\left\{\begin{array}{ll}
2,&\hbox{if $s\in\{a_0,a_1\}$;}\\
0,&\hbox{if $s\in\{a_2,a_3\}$;}\\
1,&\hbox{otherwise.}
\end{array}\right.\]
Since $\deg(f(t)^d)=3d<p-1$ for $d\in\{0,1,2,3\}$, Lemma~\ref{moments} gives
\[0=\sum_{t\in F}f(t)^d=a_0^d+a_1^d-a_2^d-a_3^d.\]
This gives a nontrivial linear dependence between four columns of the nonsingular Vandermonde matrix
$\bigl[a_i^j\,:\,0\leqslant i,j\leqslant 3\bigr]$, a contradiction.\qed
\end{pf}

\begin{corollary}\label{squares}
Let $F=\F_p$ be a field of prime order $p\geqslant5$, and suppose $\varepsilon_g=-\varepsilon_f$ where
$f(t)=t^3+ct\in F[t]$ and $g(t)=t^3+\tilde{c}t\in F[t]$. If $\tilde{c}\neq c$ then we must have $p=5$
and $\{c,\tilde{c}\}=\{2,3\}$.
\end{corollary}

\begin{pf}
If $\varepsilon_f=-\varepsilon_g$ then
\[0=\varepsilon_f+\varepsilon_g=\sum_{s\in F}(|f\inv(s)|+|g\inv(s)|)\zeta^s\]
and arguing as in the proof of Lemma~\ref{fibres}, we must have $|f\inv(s)|+|g\inv(s)|=2$ for all $s\in F$.
In particular, $|f\inv(s)|\leqslant2$ for all $s\in F$. By Lemma~\ref{triple}, $p\in\{5,7\}$.
For $p=7$ the only cubics of the form $f(t)=t^3+ct$ satisfying $|f\inv(s)|\leqslant2$ for all $s\in F$,
have $|f\inv(s)|=1,0,2,1,1,2,0$ for $s=0,1,2,3,4,5,6$ respectively; and no pair of such cubics can satisfy
$|f\inv(s)|+|g\inv(s)|=2$ for all $s\in F$.

This leaves only the case $p=5$ and the pair of cubics $f(t)=t^3+2t$, $g(t)=t^3+3t$ where
$|f\inv(s)|=1,0,2,2,0$ and $|g\inv(s)|=1,2,0,0,2$ for $s=0,1,2,3,4$ respectively.\qed
\end{pf}

\begin{theorem}\label{theorem2}
Let $p$ be an odd prime, and let $A$ be the adjacency matrix of $\Gamma(4,p)$, with
characteristic polynomial $\phi(x)=\det(xI_{p^4}-A)$.
\begin{itemize}
\item[(i)]For $p=3$, we have $\phi(x)=x^{18}(x-6)(x-3)^{12}(x+3)^{14}(x^3-9x-9)^{12}$.
\item[(ii)]For $p=5$, we have $\phi(x)=x^{220}(x-20)(x-5)^{80}(x+5)^{164}(x^2-5x-25)^{80}$.
\item[(iii)]For $5<p\equiv2\mod3$, we have
\begin{align*}\phi(x)&{}=x^{3p(p-1)}\bigl(x-p(p-1)\bigr)(x-p)^{p(p-1)^2}(x+p)^{(p-1)(2p^2-2p+1)}\\
&\quad{}\times\prod_{c=1}^{p-1}(x+p-\varepsilon_{t^3+ct}^2)^{p(p-1)^2}\end{align*}
with $p+3$ distinct roots and multiplicities as indicated by the exponents.
\item[(iv)]For $p\equiv1\mod3$, we have
\begin{align*}\phi(x)&{}=x^{3p(p-1)}\bigl(x-p(p-1)\bigr)(x-p)^{p(p-1)^2}(x+p)^{(p-1)(p^2-p+1)}\\
&\quad{}\times[(x+p-\omega_{t^3}^2)(x+p-\varepsilon_{\omega t^3}^2)(x+p-\varepsilon_{\omega^2 t^3}^2)]^{p(p-1)^2/3}\\
&{}\quad\times\prod_{a=1}^{p-1}(x+p-\varepsilon_{at^3+t}^2)^{p(p-1)^2},\end{align*}
with $p+6$ distinct roots and multiplicities as indicated by the exponents.
\end{itemize}
\end{theorem}

\begin{pf}For $p=3$ we take $R=\Z/9\Z$, $\xi=e^{2\pi i/9}$ and $\T=\{0,1,8\}$ in the notation of Section~\ref{main_proof},
and
\[\phi(x)=x^{18}(x-6)(x-3)^{12}(x+3)^{14}\prod_{c\in\T}(x+3-\varepsilon_{t^3+3ct}^2)^{12}.\]
We compute
\[\varepsilon_{t^3}=1+\xi+\xi^8,\quad\varepsilon_{t^3+3t}=1+\xi^4+\xi^5,\quad
\varepsilon_{t^3+6t}=1+\xi^2+\xi^7\]
and
\[\prod_{c\in\T}(x+3-\varepsilon_{t^3+3ct}^2)=x^3-9x-9,\]
so (i) follows.

Conclusion (iii) follows immediately from the previous results; and for $p=5$, the same
reasoning yields
\[\phi(x)=x^{60}(x-20)(x-5)^{80}(x+5)^{164}\prod_{c=1}^4(x+5-\varepsilon_{t^3+ct}^2)^{80}\]
where we need only to check for coincidence of roots. Straightforward computations show that
$\varepsilon_{t^3+t}=3+\zeta^2+\zeta^3=\frac{5-\sqrt{5}}2$ and
$\varepsilon_{t^3+4t}=2-\zeta^2-\zeta^3=\frac{5-\sqrt{5}}2$, whence
\[(x+5-\varepsilon_{t^3+t}^2)(x+5-\varepsilon_{t^3+4t}^2)=x^2-5x-25;\]
similarly, $\varepsilon_{t^3\pm2t}=\pm(1+2\zeta^2+2\zeta^3)=\mp\sqrt{5}$, giving
\[(x+5-\varepsilon_{t^3+2t}^2)(x+5-\varepsilon_{t^3+3t}^2)=x^2.\]
This yields (ii), and conclusion (iv) follows similarly from the previous results.\break\qed
\end{pf}
As an example, the spectrum of $D(4,13)$ contains $\pm\varepsilon$ where
$\varepsilon=\varepsilon_{4t^3}=1+6\cos\frac{8\pi}{13}+6\cos\frac{12\pi}{13}\approx-6.9533$.
Compare $|\varepsilon|$ with
$2\sqrt{12}\approx6.9282$ and $2\sqrt{13}\approx7.2111$ to see that while $D(4,13)$ is an 
expander, it is not quite Ramanujan. Similar conclusions are found for other values including $p=19,37$.

\section*{Acknowledgements}

This research was supported in part by NSF grant DMS-1400281. The second author is grateful to her doctoral
supervisor Dr.\ Felix Lazebnik for providing direction in this research.

\end{document}